\documentclass[amstex,12pt,amssymb]{article}
\usepackage{mathtext}
\usepackage{amsmath}
\usepackage{amssymb}
\usepackage{amsxtra}
\usepackage{latexsym}
\usepackage{ifthen}

\usepackage{indentfirst}

\usepackage{color}

\def\DD{P\ r\ o\ o\ f}

\makeatletter
\def\thmstyle{\it} 
\def\@begintheorem#1#2{\it \trivlist \item[\hskip
        \labelsep{\bf #1\ #2.}]\thmstyle}
\def\@opargbegintheorem#1#2#3{\it \trivlist \item[\hskip
        \labelsep{\bf #1\ #2\ (#3).}]\thmstyle}
\makeatother

\begin{document}

\renewcommand{\baselinestretch}{1.25}
\newcommand {\beq}{\begin{equation}}
\newcommand {\eeq}{\end{equation}}

\newtheorem{theorem}{Theorem}[section]
\newtheorem{lemma}[theorem]{Lemma}
\newtheorem{proposition}[theorem]{Proposition}
\newtheorem{corollary}[theorem]{Corollary}
\newtheorem{conjecture}[theorem]{Conjecture}
\newtheorem{definition}[theorem]{Definition}
\newtheorem{remark}[theorem]{Remark}
\newtheorem{claim}[theorem]{Claim}

\title{The local principle of large deviations for compound Poisson process  with catastrophes
 \author{\bf A. Logachov$^{1,2,3}$, O. Logachova$^{3,4}$,  A. Yambartsev$^{5}$
 }
}\maketitle

{\footnotesize
\noindent $^1$ Laboratory of Probability Theory and Mathematical Statistics, Sobolev Institute of Mathematics,
Siberian Branch of the RAS, Koptyuga str. 4, Novosibirsk 630090, RF\\
E-mail:\ \ avlogachov@mail.ru

\noindent $^2$ Novosibirsk State University,  Pirogova str. 1, Novosibirsk 630090, RF

\noindent $^3$ Novosibirsk State University of Economics and Management, Kamenskaya str. 56, Novosibirsk 630099, RF

\noindent $^4$ Siberian State University of Geo-systems and Technologies, Plakhotnogo str.
10, 630108 Novosibirsk, RF\\
E-mail:\ \ omboldovskaya@mail.ru

\noindent $^5$ Institute of Mathematics and Statistics, University of S\~ao Paulo,
1010 Rua do Mat\~ao, CEP 05508--090, S\~ao Paulo SP, Brazil \\
E-mail:\ \ yambar@ime.usp.br}

\begin{abstract} The continuous time Markov process considered in this paper belongs to a class of population models with linear growth and catastrophes. There, the catastrophes happen at the arrival times of a Poisson process, and at each catastrophe time, a randomly selected portion of the population is eliminated. For this population process, we derive an asymptotic upper bound for the maximum value and prove the local large deviation principle.
\end{abstract}



\section{Introduction}

\noindent
{\bf Motivation and historical remarks.}  First, we recall that the so-called \textit{random processes with resettings} has recently
reappeared in the context of random search models, demographic models, biological and chemical models.
See, for example, \cite{Bhat} -- \cite {Mon} and the references therein.
Informally, a {\it random process with resettings} is constructed as follows.
We modify a Markov process (i.e., a random walk, a birth-and-death process, a diffusion, etc.) by introducing the jumps (also called the {\it reset points} or {\it resettings}) to a fixed state (the origin). The jumps
happen at the arrival times of a renewal process, independent from the original Markov process.
In the above-cited papers the time intervals between resettings can be
deterministic, or they can have an exponential distribution (i.e., a Poisson process), Weibull distribution
or the distribution which depends on the current state of the process.
The cited works investigate the questions concerning the stationary distribution of random processes with resettings, its limiting behavior, and various properties of its trajectories.



Importantly, random processes with resettings are a subclass of \textit{random processes with catastrophes}
that emerged in the 1970's and 80's (see, for example, \cite{DiCGNR} -- \cite{Brock}) and stayed in the focus of researchers ever since.
A \textit{random process with catastrophes} is constructed as follows.
We modify a given Markov process by introducing additional jumps. These jumps happen
at the jump times, also called the {\it catastrophe instances}, governed by an independent renewal process.
At each jump time, a new state is selected at random, and the process restarts from that state.
In many population processes with resettings, the reset instances can be viewed
as the times of a \textit{total catastrophe}, i.e., a catastrophe which removes
instantaneously the whole population as the process jumps to the origin.


Due to a peculiarity of its applications in population dynamics and risk models,
the one dimensional random processes with catastrophes satisfy the following conditions:
(i) the processes assume only nonnegative values;
(ii) the time intervals between the catastrophes are i.i.d. and exponentially distributed;
(iii) the value of the process decreases at the catastrophe times, i.e.,
for a catastrophe time $\tau$ of the process $\xi(t)$ we have
$\xi(\tau -) > \xi(\tau)$ almost surely.

Brockwell et al. \cite{Brock} contains, apparently, the first systematic
overview of random processes with catastrophes.
We found a detailed historical overview of
results about models with catastrophes in the recent work of Ben-Ari et al. \cite{Ben},
where the authors found the invariant measure and established the
exponential convergence for a random walk with catastrophes.

To the best of our knowledge, the term \textit{uniform catastrophe} was introduced for the first time in Brockwell et al. \cite{Brock}.
There, the authors considered the population growth models subject to three types of catastrophes:
(i) geometric catastrophes, (ii) binomial catastrophes, and (iii) uniform catastrophes. The
catastrophes occur according to an independent Poisson process,
which makes the intervals between
catastrophes exponentially distributed, and the resulting process -- a continuous time Markov process.
The authors investigated the distribution of the extinction times in a simple birth process with uniform catastrophes.

The growth processes with uniform disasters were also used in modeling of the dynamics of pest population
\cite{Kyr1, Kyr2}, where damage caused by the pest population was represented by a cost function.
As the optimality criteria, the authors used the minimization of the long-run average cost per unit time.

Recently, the processes with uniform catastrophes have found their way into an interesting application
of modeling the dynamics of microbial populations \cite{microbial}.
There, the role of so-called phase variation in the dynamics of several microbial populations
was studied.
In their model, the microbial populations were evolving in an environment subjected to occasional catastrophic events.
The probability of a catastrophic event depended on the composition of the population.
In other words, the authors considered interacting
random catastrophic processes in their work. In general, we believe that the interacting random processes with catastrophes
 would be another interesting extension to consider as a future direction. There, the large deviation techniques
 could be applied in order to find some optimal strategies for population survival.

Another area of the potential application of the processes with uniform catastrophes is blockchain queueing
 models \cite{bitcoin}. There, the arrival of customers or transactions (as in blockchain models)
is a homogeneous Poisson process.
In the service protocol, at the arrival times of another homogeneous Poisson process,
a random fraction of the queue is selected, removed from the queue, and served simultaneously.
Thus, blockchain queueing models \cite{bitcoin} can be
considered from the viewpoint of the queueing models with catastrophes,
as in \cite{C95} and \cite{DE13}.

Our initial interest in the processes with uniform catastrophes came from the modeling of the dynamics of a
 spread -- the difference between the best ask and the best bid prices of an asset in the high-frequency market.
 All orders on the market can be roughly divided into two types: market orders (they can only increase the spread)
 and limit orders (they decrease the spread). Usually, large spread changes are attributed to changes in some
characteristics of the market, for example, changes in liquidity \cite{Doy}. Our aim is to understand how large
changes in the spread occur without altering parameters of the model, i.e., when the arrival rates for the market and the limit orders are constant in time, and the distributions of increments are not time dependent.
We use large deviation techniques in the analysis.

\bigskip

\noindent
{\bf Definition of the compound Poisson process with almost-uniform catastrophes.}
Further in the paper we will suppose that all random elements are defined
on probability space $(\Omega,\mathfrak{F},\bf{P})$. Consider Markov process
 $\xi (t), t\in \mathbb{R}^+$, starting from zero, $\xi (0)=0$, with state space $\mathbb{Z}^+:=\{0\}\cup\mathbb{N}$.
 The dynamics of the process $\xi(\cdot)$ is described in the following way.

If $\xi (t)=0$, then the state of random process does not change
during a random time $\tau$ exponentially distributed with rate $\alpha>0$.
At the moment $t+\tau$ the random process jumps to a state
$r$, $r\in \mathbb{Z}^+$ with probability
$$
\mathbf{P}(\xi (t+\tau)=r)=:\mathbf{P}_r,
$$
where for all $r\in \mathbb{Z}^+$ the inequality $\mathbf{P}_r\geq 0$ holds and
$\sum_{r=0}^\infty\mathbf{P}_r=1$.

If $\xi (t)=x\in\mathbb{N}$, then the process remains at the same state
during the random time $\tau$, which has exponential distribution with parameter $\alpha>0$.
At the time $t+\tau$ the process jumps to the state $x+r$, $r\in \mathbb{Z}^+$
with probability
$$
\mathbf{P}(\xi (t+\tau)=x+r)=\frac{\lambda}{\lambda+\mu}\mathbf{P}_r,
$$
and it jumps to the state $x-d$, $1\leq d\leq x$ with probability
$$
\mathbf{P}(\xi (t+\tau)=x-d)=\frac{\mu}{\lambda+\mu}\mathbf{Q}_d(x),
$$
where $\lambda>0$, $\mu>0$, for all $1\leq d\leq x$ the inequality
$\mathbf{Q}_d(x)\geq 0$ holds and $\sum_{d=1}^x\mathbf{Q}_d(x)=1$.

It is easy to see that if $\mu=0$ (there are no catastrophes), then the process $\xi(t)$
is integer valued monotonically non-decreasing compound Poisson processes.
Throughout the paper, we will suppose that the following two conditions hold for the distributions $\mathbf{P}$ and $\mathbf{Q}$:
\begin{enumerate}
 \item[$\mathbf{A}$]: there exists $r\in \mathbb{N}$ such that $\mathbf{P}_r>0$;

\item[$\mathbf{U}$]: there exists $\Delta>1$ such that for all $k\in \mathbb{N}$, $x\in \mathbb{N}$,
$1\leq d\leq x$ the inequality hold
$$
\frac{1}{\Delta x}\leq\mathbf{Q}_d(x)\leq\frac{\Delta}{ x}.
$$
\end{enumerate}
The condition {\bf U} limits the class of the distributions for the catastrophic jumps. We will call this class {\it
almost-uniform catastrophes}.

\bigskip

\noindent
{\bf Large Deviation and Local Large Deviation Principles.}
In this work, we study the limiting properties (asymptotic upper bound for the maximum value and the local large deviation principle) of an integer-valued non-decreasing compound Poisson process with almost-uniform catastrophes catastrophes.
In our earlier work \cite{LLYv2}, we proved the large deviation principle in the phase space for the homogeneous Poisson processes with uniform catastrophes. Thus, providing the optimal trajectory for large fluctuation.

As far as we know, there are no other results concerning the large deviations for the processes with catastrophes expressed in terms of the trajectories of the processes. Of the related results, we should mention here a paper of Frank den Hollander et al. \cite{Holl1},
where an additive and an integral functionals of Brownian motion with resetting were studied, and
a general variational formula for the large deviation rate function was derived. Also,
we would like to highlight the paper of Meylahn et al. \cite{Mey}, where the Large Deviation Principle (LDP) was established
for the integral functionals of the diffusion processes with resettings.

Thus, we are interested in estimating the velocity of the growth of the maximum of the process
$\xi(t)$ and in establishing the local large deviation principle (LLDP) for the family of processes
\begin{equation}\label{fam}
\xi_T(t):=\frac{\xi(Tt)}{T}, \ t\in [0,1],
\end{equation}
where $T$ is an unbounded increasing parameter.

Trajectories of the process $\xi_T(\cdot)$ belongs
almost sure to the space $\mathbb{D}[0,1]$ of c\`adl\`ag functions, i.e. the functions that are continuous from the right, and have a limit from the left.
For $f,g \in \mathbb{D}[0,1]$ let
$$
\rho(f,g)=\sup\limits_{t\in[0,1]}|f(t)-g(t)|.
$$

Recall the definition of LLDP
\begin{definition}
 A family of random processes $\xi_T(\cdot)$ satisfies
LLDP  on the set $G\subset \mathbb{D}[0,1]$ with the rate function
$I = I(f)\,:\, \mathbb{D}[0,1] \rightarrow [0,\infty]$  and the normalizing function
$\psi(T)$ such that $\lim\limits_{T\rightarrow\infty}\psi(T) = \infty$,
 if the following equality holds for any function $f \in G$
\begin{equation}\label{opr}
\begin{aligned}
&\lim_{\varepsilon\rightarrow 0}\limsup_{T\rightarrow \infty}\frac{1}{\psi(T)}
\ln\mathbf{P}(\xi_T(\cdot)\in U_\varepsilon(f))
\\
&=\lim_{\varepsilon\rightarrow 0}\liminf_{T\rightarrow \infty}\frac{1}{\psi(T)}
\ln\mathbf{P}(\xi_T(\cdot)\in U_\varepsilon(f))=-I(f),
\end{aligned}
\end{equation}
where
$
U_\varepsilon(f):=\{g\in \mathbb{D}[0,1]: \ \rho(f,g)<\varepsilon\}.
$
\end{definition}
For the more details about the LLDP notion see \cite{BorMog1}, \cite{BorMog2}.

We will also need the definition of large deviation principle (LDP).
Denote the closure and  the interior of a set $B$ by  $[B]$ and $(B)$, respectively.

\begin{definition} \label{d1.2}
  A family of random processes $\xi_T(\cdot)$ satisfies
LDP  in the metric space $(\mathbb{D}[0,1],\rho)$ with the rate function
$I = I(f)\,:\, \mathbb{D}[0,1] \rightarrow [0,\infty]$ and the normalizing function
$\psi(T)$ such that  $\lim\limits_{T\rightarrow\infty}\psi(T) = \infty$,
if for any $c \geq 0 $ \ the set $\{ f \in \mathbb{D}[0,1]\,:\, I(f) \leq c \}$  is a compact set in metric space $(\mathbb{D}[0,1],\rho)$,
 and for any set $B \in \mathfrak{B}_{(\mathbb{D}[0,1],\rho)}$ the following inequalities hold:
$$
 \limsup_{T \rightarrow \infty} \frac{1}{\psi(T)} \ln \mathbf{P}(\xi_T(\cdot) \in B ) \leq - I([B]),
$$
$$
\liminf_{T \rightarrow \infty} \frac{1}{\psi(T)} \ln \mathbf{P}(\xi_T(\cdot) \in B )\geq -I((B)),
$$
where $I(B) = \inf\limits_{y \in B} I(y)$ for $B\in \mathfrak{B}_{(\mathbb{D}[0,1],\rho)}$, where
$\mathfrak{B}_{(\mathbb{D}[0,1],\rho)}$ is the Borel $\sigma$-algebra  constructed by open
cylindrical subsets of the space $\mathbb{D}[0,1]$, and we determine $I(\emptyset) = \infty$.
\end{definition}

 Further  (Remark~\ref{noLDP}) we will see, that the family of the processes (\ref{fam}) does not satisfy
the LDP on the space $\mathbb{D}[0,1]$ of c\`adl\`ag functions equipped with Skorohod metric.
However, we hypothesize that it is possible to establish the LDP on the  $L^1[0,1]$ space. It is one of the directions of
our future research. Although we have no LDP on $\mathbb{D}[0,1]$, it is possible to establish the local large deviation principle
for our processes (see Theorem~\ref{t2.2}): the rough exponential asymptotics for the probability that the process stays in a
neighborhood of a function which
belongs to the set of absolutely continuous functions on the interval $[0,1]$ starting at zero and
taking positive values over $(0,1]$.

The above defined processes are characterized by the irregular behavior resulting from
 ``large'' downward jumps.
Interestingly, the proof of the LLDP shows that the corresponding asymptotics is the same for
the processes with high frequency ``small'' downward jumps as for the processes with exponentially small frequency of
``large'' downward jumps.
Curiously, the binomial catastrophes considered in  \cite{Brock} have these two characteristics
simultaneously: the ``small'' downward jumps and the
exponentially rare ``large'' downward jumps. In the case of the compound Poisson processes with
binomial catastrophes we expect that the rate function will be the same as in Theorem~\ref{t2.2}, but only on the space of monotone increasing
continuous functions with finite variation. However, we suspect a different rate function
on the other spaces of non-monotone increasing continuous functions with finite variations.

\bigskip
We use the following notations in the paper:
$\mathbb{AC}_0^M[0,1]$ is the set of monotonically nondecreasing and absolutely continuous
functions on the interval $[0,1]$ that start at zero;
 $\mathbb{AC}_0^+[0,1]$ is the set of absolutely continuous functions
on the interval $[0,1]$ starting at zero and taking positive values for
$t\in(0,1]$; $\text{Var}f_{[0,a]}$ is a total variation of the function
$f$ on the interval $[0,a]$; $\overline{B}$ is the complement of the set $B$;
$\mathbf{I}(B)$ is the indicator function of the set $B$;
$\lfloor a\rfloor$ is the integer part of the number $a$.

\bigskip
The paper is organized as follows. It consists of seven sections.
In Section~2 two different representations of the process $\xi(t)$ are provided;
in Section~3 the two main results (Theorem~\ref{t2.1} and Theorem~\ref{t2.2}) are formulated;
in Section~4 we prove Theorem~\ref{t2.1} about an upper bound for the growth rate of the maximum of
the process $\xi (t)$; in Section~5 Theorem~\ref{t2.2} (LLDP) is proved; in Section~6 we study a particular case of
the process with almost-uniform catastrophes where the growth process is a homogeneous Poisson process;
Section~7 contains some technical results used in the proofs the two main results.

\section{On the representations of random process $\xi (t)$}


Throughout the paper, the following two equivalent representations for the process $\xi(t)$ will be used.

\begin{enumerate}
 \item
 It is easy to see that the random process $\xi (t)$ can be represented as follows
\beq \label{2.1}
\xi(t):=\sum\limits_{v=0}^{\nu_1(t)}\gamma_v-\sum\limits_{k=0}^{\nu_2(t)}\zeta_k(\xi(\tau_k-)),
\eeq
where $\nu_1(t)$, $\nu_2(t)$ are independent Poisson processes
with parameters $\mathbf{E}\nu_1(t)=\frac{\alpha\lambda}{\lambda+\mu} t$ and
$\mathbf{E}\nu_2(t)=\frac{\alpha\mu}{\lambda+\mu} t$; $\tau_0=0$,
where $\tau_1,\dots,\tau_k,\dots$ are jump instances of the process $\nu_2(t)$;
random variables $\gamma_v$, $\zeta_k(x)$, $v\in \mathbb{Z}^+$,
$k\in \mathbb{Z}^+$, $x\in \mathbb{Z}^+$ are mutually independent and do not depend on
$\nu_1(t)$ and $\nu_2(t)$; $\zeta_0(x)=0$, $x\in \mathbb{Z}^+$;
$\mathbf{P}(\zeta_k(x)=d)=\mathbf{Q}_d(x)$, $1\leq d\leq x$, $k\in \mathbb{N}$, $x\in \mathbb{N}$;
$\mathbf{P}(\zeta_k(0)=-r)=\mathbf{P}_r$,  $r\in \mathbb{Z}^+$;
$\gamma_0=0$;
$\mathbf{P}(\gamma_v=r)=\mathbf{P}_r$, $v\in \mathbb{N}$, $r\in \mathbb{Z}^+$.
In order to make the notation simple, we denote
\begin{equation}\label{2.1.1}
\xi^+(t) :=  \sum\limits_{v=0}^{\nu_1(t)}\gamma_v \ \mbox{ and }\
\xi^-(t) := \sum\limits_{k=0}^{\nu_2(t)}\zeta_k(\xi(\tau_k-)).
\end{equation}
Then the representation (\ref{2.1}) can be rewritten as 
\beq \label{2.1.2}
\xi(t) := \xi^+(t) - \xi^-(t).
\eeq
Note that $\xi^+(t)$ is a compound Poisson process and
the term $\xi^-(t)$ accumulates the disasters.

\item Consider the Markov chain $\eta(k)$, $k\in \mathbb{Z}^+$,
with state space $\mathbb{Z}^+$ with transition probabilities
$$
{\bf P}(\eta(k+1)=j|\eta(k)=x)=
\left\{
           \begin{array}{lcl}
                              \frac{\lambda}{\lambda+\mu}\mathbf{P}_r,~~~\mbox{if}~~~j=x+r, \ \ x\neq 0,\\
\frac{\mu}{\lambda+\mu}\mathbf{Q}_d(x),~~~\mbox{if}~~~ j=x-d, \ x\neq 0, \\
\mathbf{P}_r,~~~\mbox{if}~~~j=r, x=0,
                              \end{array}
                              \right.
$$
where $\lambda$ and $\mu$ are positive constants. Set $\eta(0)=0$.

Let $\nu(t)$, $t\in \mathbb{R}^+$ be the Poisson process with parameter
$\mathbf{E}\nu(t)=\alpha t$, which does not depend on the Markov chain
$\eta(\cdot)$. Then the random process $\xi(t)$ allows the following representation
\beq \label{2.2}
\xi(t):=\eta(\nu(t)), \ t\in \mathbb{R}^+.
\eeq
\end{enumerate}

\noindent
Further we will suppose that the random variables $\gamma_v$, $v\in \mathbb{N}$
satisfy the following Carmer's condition:
\begin{enumerate}
 \item[$\mathbf{C}_\infty$]: for any $c>0$ the inequality
$
\mathbf{E}e^{c\gamma_v}<\infty
$
hold.
\end{enumerate}
%

\section{Main results}

Here we provide the formulation of the main theorems.

\begin{theorem} \label{t2.1} Let us fix a constant $b>0$. Then for any $\varepsilon>0$
the following equality holds
$$
\mathbf{P}\bigg(\lim\limits_{T\rightarrow\infty}\sup\limits_{t\in[0,1]}\frac{\xi(Tt)}{T^b}>\varepsilon\bigg)=0.
$$
\end{theorem}

Note that from Theorem \ref{t2.1} it follows that the scaled process $\xi_T(t)$
converges to zero uniformly almost surely.

Every function $f\in \mathbb{AC}_0^+[0,1]$ can be uniquely represented
as a difference of functions
$f^+\in \mathbb{AC}_0^M[0,1]$ and $f^-\in \mathbb{AC}_0^M[0,1]$
such that
$$\text{Var}f_{[0,1]}=\text{Var}f^+_{[0,1]} +\text{Var}f^-_{[0,1]}.$$
Nondecreasing functions $f^+$ and $f^-$ are called the positive
and the negative variations of the function $f$ respectively, see \cite[Ch. 1, \S 4]{Rise}.
For a given function $f\in \mathbb{AC}_0^+[0,1]$ let $B_f$ denote the set of monotonically nondecreasing functions
$g$, such that for almost all $t\in[0,1]$ the inequality $\dot{g}(t)\geq \dot{f}^{+}(t)$
holds, where $\dot{f}$ stands for the derivative of function $f$. We will use the following notations:
$$
A(y)  :=\ln \mathbf{E}e^{y\xi^+(1)}, \ \ \ \Lambda(x):=\sup\limits_{y\in \mathbb{R}}(yx-A(y)).
$$
\begin{theorem} (LLDP) \label{t2.2} Let the conditions
$\mathbf{C}_\infty$, $\mathbf{A}$ and $\mathbf{U}$ hold.
Then the family of random processes $\xi_T(\cdot)$ satisfies
LLDP on the set $\mathbb{AC}_0^+[0,1]$ with normalized function
$\psi(T)=T$ and the rate function
$$
I(f)= \frac{\alpha\mu}{\lambda+\mu}+\inf\limits_{g\in B_f}\int_0^1\Lambda(\dot{g}(t))dt.
$$
\end{theorem}

\begin{remark}\label{noLDP}
Note that it is impossible to obtain LDP for the family
$\xi_T(\cdot)$ in the metric space $(\mathbb{D}[0,1],\rho_S)$, where $\rho_S$
is the Skorokhod metric, because the corresponding
family of measures is not exponentially tight, see \cite[Remark (a), p. 8]{DZ}.
\end{remark}

\section{Proof of Theorem~\ref{t2.1}}

First, consider the case when $T$ assumes only the integer values.
$$
\begin{aligned}
\mathbf{P}\Bigl(\sup\limits_{t\in[0,1]}\xi(Tt)>T^b\varepsilon\Bigr)\leq &
\mathbf{P}\Bigl(\sup\limits_{t\in[0,1]}\xi(Tt)>T^b\varepsilon,\nu(T)\leq3\alpha T\Bigr)
\\ & +\mathbf{P}(\nu(T)>3\alpha T):=\mathbf{P}_1(T)+\mathbf{P}_2(T).
\end{aligned}
$$
We bound $\mathbf{P}_1(T)$ from above. Using the second representation (\ref{2.2}),
by virtue of independence of the Poisson process $\nu(\cdot)$ and the Markov chain $\eta(\cdot)$
we obtain
\beq \label{3.1}
\begin{aligned}
\mathbf{P}_1(T) & =\sum\limits_{n=0}^{\lfloor 3\alpha T \rfloor}\mathbf{P}\Bigl( \sup\limits_{0\leq k \leq n}\eta(k)>T^b\varepsilon\Bigr)
\mathbf{P}(\nu(T)=n)
\\
&\leq \mathbf{P}\Bigl(\sup\limits_{0\leq k \leq \lfloor 3\alpha T \rfloor}\eta(k)>T^b\varepsilon\Bigr)
\leq\sum\limits_{k=0}^{\lfloor 3\alpha T \rfloor}\mathbf{P}(\eta(k)>T^b\varepsilon).
\end{aligned}
\eeq
It is obvious that for any constant $C_1$ and for a sufficiently large value of $T$ the following equality
holds
\beq \label{3.2}
\mathbf{P} (\eta(k)>T^b\varepsilon)=\mathbf{P} \left( \eta(k)\mathbf{I}(\eta(k)>C_1)>T^b\varepsilon \right).
\eeq
Choose $C_1=C_1(\lambda,\mu,\Delta)$, where $C_1(\lambda,\mu,\Delta)$ is a constant from the Lemma~\ref{l5.1}.
Using Lemma~\ref{l5.1} and Chebyshev's inequality we obtain
\beq \label{3.3}
\begin{aligned}
\mathbf{P}\left(\eta(k)\mathbf{I}(\eta(k)>C_1)>T^b\varepsilon \right) & \leq
\mathbf{P}\left(\eta(k)\mathbf{I}(\eta(k)>C_1)>T^{1/u}\varepsilon\right)
\\
& \leq\frac{\mathbf{E}\eta^{3u}(k)\mathbf{I}(\eta(k)>C_1)}{T^3\varepsilon^{3u}}\leq\frac{C_2}{T^3\varepsilon^{3u}},
\end{aligned}
\eeq
where $u:=\min\{d\in \mathbb{N}:1/d\leq b\}$, and $C_2=C_2(\lambda,\mu,\Delta)$ is the constant
from Lemma~\ref{l5.1}.

From (\ref{3.1}), (\ref{3.2}), (\ref{3.3}) it follows that for sufficiently large $T$
\beq \label{3.4}
\mathbf{P}_1(T)\leq \frac{3\alpha C_2}{T^2\varepsilon^{3u}}.
\eeq
Next, we provide an upper bound for $\mathbf{P}_2(T)$.
Applying the Stirling's formula  it follows that for sufficiently large
$T$
\begin{equation}\label{3.5}
\begin{aligned}
\mathbf{P}_2(T)& =e^{-\alpha T}\sum\limits_{n=\lfloor 3\alpha T \rfloor+1}^\infty\frac{(\alpha T)^n}{n!}\leq
e^{-\alpha T}\sum\limits_{n=\lfloor 3\alpha T \rfloor+1}^\infty\frac{e^n(\alpha T)^n}{(3\alpha T)^n}
\\
&
\leq
e^{-\alpha T}\sum\limits_{n=0}^\infty\left(\frac{e}{3}\right)^n=e^{-\alpha T}\frac{3}{3-e}.
\end{aligned}
\end{equation}

From (\ref{3.4}) and (\ref{3.5}) it follows that the series
$
\sum\limits_{T=1}^\infty\mathbf{P}\Bigl(\sup\limits_{t\in[0,1]}\frac{\xi(Tt)}{T^b}>\varepsilon\Bigr)
$
converges. Therefore, by virtue of the Borel-Cantelli lemma, Theorem~\ref{t2.1}
is proved for the case when the parameter $T$ takes only integer values.

We will show now that Theorem~\ref{t2.1} holds for any real $T>0$. There, for all $T> 0$, we have
$$
0\leq\sup\limits_{t\in[0,1]}\frac{\xi(Tt)}{T^b}\leq \sup\limits_{t\in[0,1]}\frac{\xi((\lfloor T \rfloor+1)t)}{T^b}=
\frac{(\lfloor T \rfloor+1)^b}{T^b}\sup\limits_{t\in[0,1]}\frac{\xi((\lfloor T \rfloor+1)t)}{(\lfloor T \rfloor+1)^b}.
$$
Hence, due to the fact that $\sup\limits_{t\in[0,1]}\frac{\xi((\lfloor T \rfloor+1)t)}{(\lfloor T \rfloor+1)^b}$
converges a.s. to zero, Theorem~\ref{t2.1} is proved. $\Box$

\section{Proof of Theorem~\ref{t2.2}}

\noindent
We will use representation (\ref{2.1.2}) to prove the theorem.
Applying the scaling (\ref{fam}) we denote
\beq \label{4.2}
\xi_T(t) = \xi^+_T(t) - \xi^-_T(t),
\eeq
where
$$
\xi^+_T(t) := \frac{\xi^+(Tt)}{T} \ \mbox{ and }\ \xi^-_T(t) := \frac{\xi^-(Tt)}{T}.
$$

First, we bound $\mathbf{P}(\xi_T(\cdot)\in U_\varepsilon(f))$ from above.
For all $c>0$ and $\delta>0$ we have
$$
\begin{aligned}
\mathbf{P}(\xi_T(\cdot)\in U_\varepsilon(f))\leq & \mathbf{P}
\bigg(\sup\limits_{t\in[\delta,1]}|\xi_T(t)-f(t)|<\varepsilon,A_c\bigg)
\\
+ & \mathbf{P}\bigg(\sup\limits_{t\in[\delta,1]}|\xi_T(t)-f(t)|<\varepsilon,\overline{A_c}\bigg)
:=\mathbf{P}_1+\mathbf{P}_2,
\end{aligned}
$$
where
$$
A_c:=\{\omega:\nu_2(T)-\nu_2(\delta T)\leq cT\}.
$$
Now, we provide an upper bound for $\mathbf{P}_1$. Using the formula (\ref{4.2}), we obtain
$$
\mathbf{P}_1=\mathbf{P}\bigg(\sup\limits_{t\in[\delta,1]}|\xi_{T}^+(t)-\xi_{T}^-(t)-f(t)|<\varepsilon,A_c\bigg).
$$
Denote
$
m_\delta:=\min\limits_{t\in[\delta,1]}f(t).
$ Note that for sufficiently small $\varepsilon$ the inequality $m_\delta>\varepsilon$ is
satisfied. Therefore, if  $\sup\limits_{t\in[\delta,1]}|\xi_{T}(t)-f(t)|<\varepsilon$, then
$\xi_{T}(t)>0$ for $t\in[\delta,1]$, and
the random process $\xi_{T}^-(t)$ does not monotonically decrease on $[\delta,1]$.
For any $r>0$ the following inequality holds
$$
\begin{aligned}
\mathbf{P}_1 & =\mathbf{P}\Bigl(\sup\limits_{t\in[\delta,1]}|\xi_{T}^+(t)-\xi_{T}^-(t)-f(t)|<\varepsilon,A_c\Bigr)
\\
& \leq
\mathbf{P}\Bigl(\sup\limits_{t\in[\delta,1]}|\xi_{T}^+(t)-\xi_{T}^-(t)-f(t)|<\varepsilon,A_c,\xi_{T}^+(\cdot)\in K_r^\varepsilon\Bigr)
+\mathbf{P}\left(\xi_{T}^+(\cdot)\in \overline{K_r^\varepsilon}\right):=\mathbf{P}_{11}+\mathbf{P}_{12},
\end{aligned}
$$
where
$$
K_r^\varepsilon:=
\bigg\{v\in \mathbb{D}_0^M[0,1]:\inf\limits_{g\in K_r}\sup\limits_{t\in[0,1]}|g(t)-v(t)|\leq\varepsilon\bigg\},
$$
$\mathbb{D}_0^M[0,1]$  is the set of monotonically non-decreasing c\`adl\`ag functions
starting at zero,
$$
K_r :=\big\{g:I_1(g)\leq r\big\} \ \ \mbox{ and } \ \
I_1(g):=\int_0^1\Lambda(\dot{g}(t))dt,
$$
Note that, from Theorem~3.3 in \cite{BorMog3}, it follows that the set $K_r$ is a compact set.

Next, we bound  $\mathbf{P}_{11}$ from above. Denote
$$
B_f:=\big\{g\in \mathbb{D}^M_0[0,1]: \dot{g}(t)\geq \dot{f}^+(t) \text{ for almost all } t\in[0,1]\big\}.
$$

Since the random process $\xi_{T}^-(t)$ is nondecreasing on the interval $[\delta,1]$,
and $K_r$ is a compact set, Lemma~\ref{l5.21} implies that there exists a positive function
$\gamma(\cdot)>0$ such that $\gamma(\varepsilon)\rightarrow 0$ as $\varepsilon\rightarrow 0$ and
$$
\mathbf{P}_{11} \leq \mathbf{P} \left( \xi_{T}^+(\cdot)\in B_f^{\delta,\gamma(\varepsilon)},A_c,K_ r^\varepsilon \right)\leq
\mathbf{P}\left( \xi_{T}^+(\cdot)\in B_f^{\delta,\gamma(\varepsilon)},A_c \right),
$$
where
$$
B_f^{\delta,\gamma(\varepsilon)}:=\Bigl\{v\in \mathbb{D}[0,1]:\inf\limits_{g\in B_f}\sup\limits_{t\in[\delta,1]}|g(t)-v(t)|\leq\gamma(\varepsilon)\Bigr\}.
$$

Since the random processes $\xi_{T}^+(t)$ and $\xi_{T}^-(t)$ are independent, we obtain
$$
\mathbf{P}_{11}\leq\mathbf{P} \left(\xi_{T}^+(\cdot)\in B_f^{\delta,\gamma(\varepsilon)},A_c \right)=\mathbf{P} \left(\xi_{T}^+(\cdot)\in B_f^{\delta,\gamma(\varepsilon)} \right)\mathbf{P}(A_c).
$$

It means that for all $r>0$,
\beq \label{4.3}
\mathbf{P}_{1}\leq
\mathbf{P} \left( \xi_{T}^+(\cdot)\in B_f^{\delta,\gamma(\varepsilon)} \right)\mathbf{P}(A_c)+\mathbf{P}_{12}
=\mathbf{P} \left(\xi_{T}^+(\cdot)\in B_f^{\delta,\gamma(\varepsilon)} \right)\mathbf{P}(A_c)+\mathbf{P} \left(\xi_{T}^+(\cdot)\in\overline{K_ r^\varepsilon} \ \right).
\eeq

Next we provide an upper bound for $\mathbf{P}_2$.
Denote $\tau_{k_1},\dots,\tau_{k_{\lfloor cT \rfloor}}$ the first $\lfloor cT \rfloor$ jumps of the process
$\nu_2(Tt)$ belonging to the interval $[\delta,1]$. Denote
$$
\begin{aligned}
G_{k_l} & :=\{\omega:\xi(\tau_{k_l}-)\in[T(f(\tau_{k_l})-\varepsilon);T(f(\tau_{k_l})+\varepsilon)]\}, \ 1\leq l \leq \lfloor cT \rfloor,
\\
H_{k_l} &:=\{\omega:\zeta_{k_l}(\xi(\tau_{k_l}-))< 2 T\varepsilon\}, \ 1\leq l \leq \lfloor cT \rfloor.
\end{aligned}
$$
If the trajectory of the process $\xi_T(t)$ does not leave the set
$U_\varepsilon(f)$, then
$\zeta_{k_l}(\xi(\tau_{k_l}-))< 2 T\varepsilon$ for $\tau_{k_l}\in[\delta,1]$,
$1\leq l \leq \lfloor cT \rfloor$.
Therefore, the inequality
$$
\begin{aligned}
\mathbf{P}_2 &=\mathbf{P}\bigg(\sup\limits_{t\in[\delta,1]}|\xi_T(t)-f(t)|<\varepsilon,\overline{A_c}\bigg)
\\
& \leq\sum\limits_{r=\lfloor cT \rfloor}^{\infty}\mathbf{P}\bigg(\bigcap\limits_{l=1}^{\lfloor cT \rfloor}H_{k_l},\bigcap\limits_{l=1}^{\lfloor cT \rfloor}G_{k_l} \
\bigg| \ \nu_2(T)- \nu_2(\delta T)=r\bigg)
\mathbf{P}\bigl(\nu_2(T)- \nu_2(\delta T)=r\bigr).
\end{aligned}
$$
Using Lemma~\ref{l5.6}, we obtain
$$
\mathbf{P}_2\leq
\sum\limits_{r=\lfloor cT \rfloor}^{\infty}\bigg(\frac{\lfloor 2 T\varepsilon\rfloor }{\lfloor T(m_\delta-\varepsilon) \rfloor}\bigg)^{\lfloor cT \rfloor}
\mathbf{P}\bigl(\nu_2(T)- \nu_2(\delta T)=r\bigr)
\leq\bigg(\frac{\Delta\lfloor 2T\varepsilon \rfloor}{\lfloor T(m_\delta-\varepsilon) \rfloor}\bigg)^{\lfloor cT \rfloor}
$$
For sufficiently small $\varepsilon$ the inequality
$m_\delta>\sqrt{\varepsilon}$ is satisfied, hence for sufficiently large $T$
\beq \label{4.4}
\mathbf{P}_2
\leq\bigg(\frac{\Delta\lfloor 2T\varepsilon \rfloor}{\lfloor T(m_\delta-\varepsilon) \rfloor}\bigg)^{\lfloor cT \rfloor}\leq
\bigg(\frac{4\Delta\sqrt{\varepsilon}}{1-\sqrt{\varepsilon}}\bigg)^{\lfloor cT \rfloor}.
\eeq
From inequality (\ref{4.4}) it follows that for any $c>0$,
\beq \label{4.41}
\lim\limits_{\varepsilon\rightarrow 0}\limsup\limits_{T\rightarrow \infty}\frac{1}{T}\ln\mathbf{P}_2
\leq c \lim\limits_{\varepsilon\rightarrow 0}\ln\bigg(\frac{4\Delta\sqrt{\varepsilon}}{1-\sqrt{\varepsilon}}\bigg)=
-\infty.
\eeq
From Theorem 3.3 \cite{BorMog3} it follows that for any $\varepsilon>0$,
\beq \label{new3}
\lim\limits_{r\rightarrow \infty}\limsup\limits_{T\rightarrow \infty}\frac{1}{T}\ln\mathbf{P} \left(\xi_{T}^+(\cdot)\in\overline{K_ r^\varepsilon} \right)=
-\infty.
\eeq
Therefore, using (\ref{4.3}), (\ref{4.41}), (\ref{new3}), Theorem 3.3 \cite{BorMog3},
Lemma \ref{l5.4}
and the fact that the set $B_f^{\delta,\gamma(\varepsilon)}$  is closed, for all $c\in(0,1)$ and $\delta>0$
we have
$$
\begin{aligned}
& \lim\limits_{\varepsilon\rightarrow 0}\limsup\limits_{T\rightarrow \infty}\frac{1}{T}  \ln\mathbf{P}(\xi_T(\cdot)\in U_\varepsilon(f))\leq
\lim\limits_{\varepsilon\rightarrow 0}\limsup\limits_{T\rightarrow \infty}\frac{1}{T}\ln(\mathbf{P}_{11}+\mathbf{P}_{12}+\mathbf{P}_2)
\\
&\le
\lim\limits_{\varepsilon\rightarrow 0}\limsup\limits_{T\rightarrow \infty}\frac{1}{T}\ln(3\max\{\mathbf{P}_{11},\mathbf{P}_{12},\mathbf{P}_2\})
\\ & \leq
\lim\limits_{\varepsilon\rightarrow 0}\bigg(-I_1(B_f^{\delta,\gamma(\varepsilon)})
-\frac{\alpha\mu(1-\delta)}{\lambda+\mu}
 +\frac{\alpha\mu(1-\delta)c}{\lambda+\mu} - c\ln c\bigg)
\\
&=-I_1(B_f^{\delta})
-\frac{\alpha\mu(1-\delta)}{\lambda+\mu}
 +\frac{\alpha\mu(1-\delta)c}{\lambda+\mu} - c\ln c,
 \end{aligned}
$$
where
$$
B_f^{\delta}:=\bigg\{v\in \mathbb{D}[0,1]:\inf\limits_{g\in B_f}\sup\limits_{t\in[\delta,1]}|g(t)-v(t)|=0\bigg\}.
$$
Taking the limits $\delta\rightarrow 0$ and $c\rightarrow 0$ we obtain
$$
\lim\limits_{\varepsilon\rightarrow 0}\limsup\limits_{T\rightarrow \infty}
\frac{1}{T}\ln\mathbf{P}(\xi_T(\cdot)\in U_\varepsilon(f))\leq-I_1(B_f)-\frac{\alpha\mu}{\lambda+\mu}.
$$

Now we bound  $\mathbf{P}(\xi_T(\cdot)\in U_\varepsilon(f))$ from below.
$$
\mathbf{P}_3:=\mathbf{P}\bigg(\sup\limits_{t\in[0,1]}|\xi_{T}^+(t)-\xi_{T}^-(t)-f(t)|<\varepsilon\bigg)
\geq\mathbf{P}\bigg(\xi_{T}^+(\cdot)\in U_{\frac{\varepsilon}{2}}(g^*),
\xi_{T}^-(\cdot)\in U_{\frac{\varepsilon}{2}}(g^*-f)\bigg),
$$
where $g^*$ is such function that $\inf\limits_{g\in B_f}I_1(g)=I_1(g^*)$.
Note that this function exists because $$\inf\limits_{g\in B_f}I_1(g)=
\inf\limits_{g\in B_f\cap\{g:I_1(g)\leq I_1(f^+)\}}I_1(g),$$
and the set  $B_f\cap\{g:I_1(g)\leq I_1(f^+)\}$ is a compact set according Definition~\ref{d1.2}.

Since the function $g^*\in B_f$, the function $g^*-f\in \mathbb{AC}_0^M[0,1]$.
If $g^*-f\equiv 0$, then
$$
\mathbf{P}\left(\xi_{T}^+(\cdot)\in U_{\frac{\varepsilon}{2}}(g^*),
\xi_{T}^-(\cdot)\in U_{\frac{\varepsilon}{2}}(g^*-f)\right)\geq
\mathbf{P}\left(\xi_{T}^+(\cdot)\in U_{\frac{\varepsilon}{2}}(g^*),
\nu_2(T)=0\right).
$$
Therefore, since the random processes $\xi_{T}^+(\cdot)$ and $\nu_2(\cdot)$ are independent, we obtain
\beq \label{4.5}
\mathbf{P}_3\geq \mathbf{P}\left(\xi_{T}^+(\cdot)\in U_{\frac{\varepsilon}{2}}(g^*)\right)e^{-\frac{\alpha\mu}{\lambda+\mu} T}.
\eeq
Let $g^*-f\not\equiv 0$. Define
$$n(\varepsilon):=\min\bigg\{n\in \mathbb{N}:\frac{M}{n}\leq\frac{\varepsilon}{8}\bigg\},$$
where $M:=\max\limits_{t\in[0,1]}(g^*(t)-f(t))=g^*(1)-f(1)$.

Since the function $g^*-f$ is continuous and monotonically nondecreasing, there
exists a finite set of points $0=t_0<t_1<\dots<t_{n(\varepsilon)}=1$, such that the equalities
$$
g^*(t_1)-f(t_1)=\frac{M}{n(\varepsilon)}, \ g^*(t_1)-f(t_1)=\frac{2M}{n(\varepsilon)},
\dots, \ g^*(t_{n(\varepsilon)})-f(t_{n(\varepsilon)})=M
$$
are satisfied. Therefore, if the random process $\nu_2(Tt)$ does not have jumps within the interval $[0,t_1]$
and has only one jump on each of the segments $[t_{k-1},t_{k}]$, $2 \leq k \leq n(\varepsilon)$,
and if the random variables $\zeta_k(\xi(\tau_k-))$ are all in the interval
$$
\bigg(\frac{TM}{n(\varepsilon)}-2T\varepsilon^3;\frac{TM}{n(\varepsilon)}-T\varepsilon^3\bigg),
$$
then for sufficiently small $\varepsilon$  the inequality
$$
\sup\limits_{t\in[0,1]}
\bigg|\frac{1}{T}\sum\limits_{k=0}^{\nu_2(Tt)}\zeta_k(\xi(\tau_k-))-(g^*(t)-f(t))\bigg|<\frac{\varepsilon}{2}.
$$
holds. Hence, for sufficiently small $\varepsilon>0$ the following inequality
\small\beq \label{4.6}
\mathbf{P}_3\geq\mathbf{P}\bigg(\xi_{T}^+(\cdot)\in U_{\varepsilon^3}(g^*),
\xi_{T}^-(\cdot)\in U_{\frac{\varepsilon}{2}}(g^*-f)\bigg)
\geq\mathbf{P}\bigg(\xi_{T}^+(\cdot)\in U_{\varepsilon^3}(g^*),
\bigcap\limits_{k=1}^{n(\varepsilon)} A_k,\bigcap\limits_{k=1}^{n(\varepsilon)-1} B_k\bigg),
\eeq\normalsize
holds, where
$$
\begin{aligned}
A_1 & :=\{\omega:\nu_2(Tt_1)=0\}, \ A_k:=\{\omega:\nu_2(Tt_k)-\nu_2(Tt_{k-1})=1\},
\ 2 \leq k \leq n(\varepsilon),
\\
B_k & :=\bigg\{\omega:\zeta_k(\xi(\tau_k-))\in
\bigg(\frac{TM}{n(\varepsilon)}-2T\varepsilon^3;\frac{TM}{n(\varepsilon)}-T\varepsilon^3\bigg)\bigg\},
 \ 1 \leq k \leq n(\varepsilon)-1.
\end{aligned}
$$
From the inequality (\ref{4.6}) it follows that
$$
\mathbf{P}_3\geq\mathbf{P}\bigg(\bigcap\limits_{k=1}^{n(\varepsilon)-1} B_k \  \bigg| \
\xi_{T}^+(\cdot)\in U_{\varepsilon^3}(g^*),
\bigcap\limits_{k=1}^{n(\varepsilon)} A_k\bigg)
\mathbf{P}\bigg(\xi_{T}^+(\cdot)\in U_{\varepsilon^3}(g^*),
\bigcap\limits_{k=1}^{n(\varepsilon)} A_k\bigg).
$$
And from Lemma~\ref{5.5} it follows that
$$
\mathbf{P}_3\geq\bigg(\frac{\varepsilon^3}{2\Delta(g^*(1)+\varepsilon^3)}\bigg)^{n(\varepsilon)-1}\mathbf{P}\bigg(\xi_{T}^+(\cdot)\in U_{\varepsilon^3}(g^*),
\bigcap\limits_{k=1}^{n(\varepsilon)} A_k\bigg).
$$
Since the random processes $\xi_{T}^+(t)$ and $\nu_2(Tt)$ are independent, we have
\beq \label{4.7}
\begin{aligned}
& \mathbf{P}_3\geq\bigg(\frac{\varepsilon^3}{2\Delta(g^*(1)+\varepsilon^3)}\bigg)^{n(\varepsilon)-1}\mathbf{P}(\xi_{T}^+(\cdot)\in U_{\varepsilon^3}(g^*))\mathbf{P}\bigg(
\bigcap\limits_{k=1}^{n(\varepsilon)} A_k\bigg)
\\
& =\bigg(\frac{\varepsilon^3}{2\Delta(g^*(1)+\varepsilon^3)}\bigg)^{n(\varepsilon)-1}\mathbf{P}(\xi_{T}^+(\cdot)\in U_{\varepsilon^3}(g^*))
\bigg(\frac{\alpha\mu}{\lambda+\mu}T\bigg)^{n(\varepsilon)-1}e^{-\frac{\alpha\mu}{\lambda+\mu}T}
\prod\limits_{k=2}^{n(\varepsilon)}(t_k-t_{k-1}).
\end{aligned}
\eeq
Using the inequalities (\ref{4.5}),  (\ref{4.7}) and Theorem 3.3 \cite{BorMog3}, we obtain
$$
\liminf_{T\rightarrow \infty}\frac{1}{T}
\ln\mathbf{P}(\xi_T(\cdot)\in U_\varepsilon(f))\geq -\frac{\alpha\mu}{\lambda+\mu}-I_1(U_{\varepsilon^3}(g^*)).
$$
From Theorem 3.1 \cite{BorMog3} it follows that
$$
\lim\limits_{\varepsilon\rightarrow 0}\liminf_{T\rightarrow \infty}\frac{1}{T}
\ln\mathbf{P}(\xi_T(\cdot)\in U_\varepsilon(f))\geq \lim\limits_{\varepsilon\rightarrow 0}\bigg(-\frac{\alpha\mu}{\lambda+\mu}-I_1(U_{\varepsilon^3}(g^*))\bigg)=
-\frac{\alpha\mu}{\lambda+\mu}-I_1(g^*).
$$
$\Box$

\section{Model Example}

In this section we consider a common Poisson process with uniform catastrophes.
Considera discrete time Markov chain $\eta(k)$, $k\in \mathbb{Z}^+$, $\mathbb{Z}^+=\{0\}\cup\mathbb{N}$,
with state space $\mathbb{Z}^+$ and with transition probabilities
\beq
{\bf P}(\eta(k+1)=j|\eta(k)=i)=
\left\{
           \begin{array}{lcl}
                              \frac{\lambda}{\lambda+\mu},~~~\mbox{if}~~~j=i+1,\\
\frac{\mu}{i(\lambda+\mu)},~~~\mbox{если}~~~0\leq j<i, i\neq 0,\\
1,~~~\mbox{if}~~~j=1, i=0,
                              \end{array}
                              \right.
\eeq
where $\lambda$ and $\mu$ are positive constant. Let $\eta(0)=0$.

Consider a Poisson process $\nu(t)$, $t\in \mathbb{R}^+$, with parameter
$\mathbf{E}\nu(t)=\alpha t$, that does not depend on the Markov chain
$\eta(\cdot)$. Define the random process
$$
\tilde{\xi}(t):=\eta(\nu(t)), \ t\in \mathbb{R}^+,
$$
and  denote
$$
\tilde{\xi}_T(t):=\frac{\tilde{\xi}(Tt)}{T}.
$$

\begin{theorem} (LLDP) \label{t6.2}  The family of random processes $\tilde{\xi}_T(\cdot)$ satisfies
LLDP on the set $\mathbb{AC}_0^+[0,1]$ with the normalized function
$\psi(T)=T$ and the rate function
$$
I(f)= \frac{\alpha\mu}{\lambda+\mu}+\inf\limits_{g\in B_f}\int_0^1\left(\dot{g}(t)\ln\left(\frac{\dot{g}(t)(\lambda+\mu)}{\alpha\lambda}\right)-\dot{g}(t)+
\frac{\alpha\lambda}{\lambda+\mu}\right)dt.
$$
\end{theorem}
\noindent
\DD.
The proof follows from Theorem~\ref{t2.2} and Lemma \ref{l5.3}.
$\Box$

\section{Auxiliary results}

Here we will prove several auxiliary technical results.

Introduce the random variable $\gamma$ which distribution coincides with distribution of random variable
$\gamma_v$, $v\in \mathbb{N}$. Denote
$$
\begin{aligned}
k_1& 
:= \left(\mathbf{E}\gamma^{3u}\right)^{1/3u},
\\
k_2&:=\left(\left(\frac{4\Delta(\lambda+\mu)}{4\lambda\Delta+\mu (4\Delta-2)}\right)^{1/3u}-1\right)^{-1}, \ \ \ k_3:= \left(\frac{4\lambda\Delta+\mu(4\Delta-2)}{4\lambda\Delta+\mu (4\Delta-3)}\right)^{1/3u},
\\
C_1&:=C_1(\lambda,\mu,\Delta)=\left\lfloor k_1\max (k_2,k_3)\right\rfloor+1.
\end{aligned}
$$

\begin{lemma} \label{l5.1} Let $u\in \mathbb{N}$ be a fixed constant.
Then for all $k\in \mathbb{Z}^+$, the inequality
\beq \label{5.1}
 \mathbf{E}\eta^{3u}(k)\mathbf{I}(\eta(k)>C_1)\leq C_2
\eeq
holds, where
$$
C_2:=C_2(\lambda,\mu,\Delta)=\frac{4\lambda\Delta+\mu (4\Delta-3)}{\mu}(C_1)^{3u}.
$$
\end{lemma}
\noindent
\DD. The proof will be carried out by the method of mathematical induction. If $k=0$, then it
is obvious that the inequality
(\ref{5.1}) holds. Suppose that the inequality (\ref{5.1}) holds for $k=m-1$. We will show that, it holds for $k=m$.
We have
\beq \label{5.2}
\begin{aligned}
 &\mathbf{E}\eta^{3u}(m)\mathbf{I}(\eta(m)> C_1)=
 \mathbf{E} \bigl( \mathbf{E}(\eta^{3u}(m)\mathbf{I}(\eta(m)> C_1) \mid \eta(m-1)) \bigr)
\\
&=\mathbf{E}\sum\limits_{r=0}^{\infty}\mathbf{I}(\eta(m-1)=r)\mathbf{E}(\eta^{3u}(m)\mathbf{I}(\eta(m)> C_1) \mid \eta(m-1)=r).
\end{aligned}
\eeq
We will show now that
\beq \label{5.3}
\mathbf{E}(\eta^{3u}(m)\mathbf{I}(\eta(m)> C_1) \mid \eta(m-1)=r)\leq (\max(C_1,r))^{3u}\frac{4\lambda\Delta+\mu(4\Delta-3)}{4\lambda\Delta+\mu (4\Delta-2)}.
\eeq
If $r=0$ then we have
$$
\mathbf{E}(\eta^{3u}(m)\mathbf{I}(\eta(m)> C_1) \mid \eta(m-1)=0)\leq\mathbf{E}\gamma^{3u}\leq k_1^{3u}\leq C_1^{3u}\frac{4\lambda\Delta+\mu(4\Delta-3)}{4\lambda\Delta+\mu (4\Delta-2)}.
$$
If $r>0$ then from the definition of Markov chain $\eta(\cdot)$ we have
\begin{eqnarray}\label{aux1}
&& \mathbf{E}(\eta^{3u}(m)\mathbf{I}(\eta(m)> C_1) \mid \eta(m-1)=r)\leq
\mathbf{E}(\eta^{3u}(m) \mid \eta(m-1)=r) \nonumber
\\
&& \ \ \ =\frac{\lambda}{\lambda+\mu}\mathbf{E}(r+\gamma)^{3u}+\frac{\mu}{(\lambda+\mu)}\sum\limits_{d=1}^{r}(r-d)^{3u}\mathbf{Q}_d(r) 
\end{eqnarray}
Note that by Jensen inequality the following $(\mathbf{E}(\gamma^d))^{1/d} \le (\mathbf{E}(\gamma^{3u}))^{1/3u}$
holds for any $1\le d \le 3u$. Thus, the first term from (\ref{aux1}) can be bounded as
$$
\mathbf{E}(r+\gamma)^{3u} \le (1+ k_1)^{3u}.
$$
The series form (\ref{aux1}) we bound from above in the following way
$$
\begin{aligned}
&\sum\limits_{d=1}^{r}(r-d)^{3u}\mathbf{Q}_d(r)  = \sum\limits_{d=1}^{r}(r-d)^{3u}\left(\mathbf{Q}_d(r)-\frac{1}{\Delta r}\right)+ \frac{1}{\Delta r}\sum\limits_{d=1}^{r}(r-d)^{3u}
\\
&  \le (r-1)^{3u} \left(1-\frac{1}{\Delta}\right) + \frac{(r-1)^{3u-3}}{\Delta r} \sum\limits_{d=1}^{r}(r-d)^{3} 
\\
& = (r-1)^{3u} \left(1-\frac{1}{\Delta}\right) + \frac{r(r-1)^{3u-1}}{4\Delta} \le 
\frac{(r+k_1)^{3u}(4\Delta-3)}{4\Delta},
\end{aligned}
$$
where we used the formula for sum of cubes of natural numbers. Thus,
$$
\begin{aligned}
& \mathbf{E}(\eta^{3u}(m)\mathbf{I}(\eta(m)> C_1) \mid \eta(m-1)=r) \leq 
\frac{\lambda}{\lambda+\mu}(r+k_1)^{3u} +\frac{\mu(r+k_1)^{3u}(4\Delta-3)}{4(\lambda+\mu)\Delta}
\\
& \ \ \ =
\frac{4\lambda\Delta+\mu(4\Delta-3)}{4(\lambda+\mu)\Delta}(r+k_1)^{3u} 
\leq\frac{4\lambda\Delta+\mu(4\Delta-3)}{4\lambda\Delta+\mu (4\Delta-2)}(\max(C_1,r))^{3u}.
\end{aligned}
$$
From (\ref{5.2}) and (\ref{5.3}) and the inductive assumption it follows that
$$
\begin{aligned}
& \mathbf{E}\eta^{3u}(m)\mathbf{I}(\eta(m)> C_1)\leq
\frac{4\lambda\Delta+\mu(4\Delta-3)}{4\lambda\Delta+\mu (4\Delta-2)}\mathbf{E}\sum\limits_{r=0}^{\infty}(\max(C_1,r))^{3u}
\mathbf{I}(\eta(m-1)=r)
\\
&=\frac{4\lambda\Delta+\mu(4\Delta-3)}{4\lambda\Delta+\mu (4\Delta-2)}\left(\mathbf{E}\sum\limits_{r=0}^{C_1}(C_1)^{3u}
\mathbf{I}(\eta(m-1)=r)+\mathbf{E}\sum\limits_{r=C_1+1}^{\infty}r^{3u}
\mathbf{I}(\eta(m-1)=r)\right)
\\
&\leq\frac{4\lambda\Delta+\mu(4\Delta-3)}{4\lambda\Delta+\mu (4\Delta-2)}\left((C_1)^{3u}
+\mathbf{E}\eta^{3u}(m-1)\mathbf{I}(\eta(m-1)> C_1)\right)
\leq \frac{4\lambda\Delta+\mu (4\Delta-3)}{\mu}(C_1)^{3u}.
\end{aligned}
$$
$\Box$

\begin{lemma} \label{l5.2}
Let the function $f\in \mathbb{AC}_0^+[0,1]$ be represented in the form
$$
f(t)=g_1(t)-g_2(t),
$$
where  $g_1\in \mathbb{D}_0^M[0,1]$ and $g_2\in \mathbb{D}_0^M[0,1]$.
Then the inequality $\dot{g}_1(t)\geq \dot{f}^+(t)$
holds for almost all $t\in[0,1]$.
\end{lemma}

\noindent
\DD.  Assume not. Then there exist $0\le t_1<t_2 \le 1$,
 such that $g_1(t_2)-g_1(t_1)<f^+(t_2)-f^+(t_1)$.
We note that in this case the inequality $g_2(t_2)-g_2(t_1)<f^-(t_2)-f^-(t_1)$
also holds.

Since the variation of the sum of two functions does not exceed
the sum of their variations,
$$
\text{Var}g_{1[t_1,t_2]}+\text{Var}g_{2[t_1,t_2]}=
g_1(t_2)-g_1(t_1)+g_2(t_2)-g_2(t_1)\geq\text{Var}f_{[t_1,t_2]}.
$$
On the other hand
$$
g_1(t_2)-g_1(t_1)+g_2(t_2)-g_2(t_1)<f^+(t_2)-f^+(t_1)+f^-(t_2)-f^-(t_1)=\text{Var}f_{[t_1,t_2]}.
$$
The obtained contradiction completes the proof. $\Box$

\bigskip
Consider a family of functions $u_T(t)$, $t\in[0,1]$, $T>0$, such that $u_T(t)$ can be represented by
$u_T(t):=\tilde{u}_T(t)-\hat{u}_T(t)$, where
$\hat{u}_T\in \mathbb{D}_0^M[0,1]$, $\tilde{u}_T\in \mathbb{D}_0^M[0,1]\cap K_{(\mathbb{D}[0,1],\rho)} $,
and $K_{(\mathbb{D}[0,1],\rho)}\subset (\mathbb{D}[0,1],\rho)$ is a compact set.

\begin{lemma} \label{l5.21}
Suppose for a function $f\in \mathbb{AC}_0^+[0,1]$ we have
\beq \label{new2}
\lim\limits_{T\rightarrow\infty}\sup\limits_{t\in[0,1]}|u_T(t)-f(t)|=0.
\eeq
Then
$$
\lim\limits_{T\rightarrow\infty}\inf_{g\in B_f}\sup\limits_{t\in[0,1]}|\tilde{u}_T(t)-g(t)|=0.
$$
\end{lemma}

\noindent
\DD. \ Proof by contradiction. Suppose not. Then for any
$M>0$, there exist $T>M$ and $\gamma>0$ such that
\beq \label{new1}
\inf_{g\in B_f}\sup\limits_{t\in[0,1]}|\tilde{u}_T(t)-g(t)|\geq \gamma.
\eeq
Since the family of functions $\tilde{u}_T$ contains in a compact set,
the inequality (\ref{new1}) implies that there exists subsequence
$T_M$ and a function $\tilde{g}$ such that
$$
\lim\limits_{M\rightarrow\infty}\sup\limits_{t\in[0,1]}|\tilde{u}_{T_M}(t)-\tilde{g}(t)|=0, \ \ \
\inf_{g\in B_f}\sup\limits_{t\in[0,1]}|\tilde{g}(t)-g(t)|\geq \gamma.
$$
Therefore, equation (\ref{new2}) implies
$$
\lim\limits_{M\rightarrow\infty}\sup\limits_{t\in[0,1]}|\hat{u}_{T_M}(t)-(\tilde{g}(t)-f(t))|=0.
$$
Since $\hat{u}_T\in \mathbb{D}_0^M[0,1]$, the function
$$
\hat{g}(t):=\tilde{g}(t)-f(t)
$$
should belong to the set $\mathbb{D}_0^M[0,1]$.

It means that $f(t)=\tilde{g}(t)-\hat{g}(t)$, where $\tilde{g}\not\in B_f$, $\hat{g}\in\mathbb{D}_0^M[0,1]$,
which contradicts Lemma~\ref{l5.2}. $\Box$

\begin{lemma} \label{l5.3}
The family of processes $\frac{\nu_1(Tt)}{T}$ satisfies LDP
on the metric space $(\mathbb{D}[0,1],\rho)$
with the normalizing function
$\psi(T)=T$ and the rate function
$$
I_1(f)= \left\{
\begin{aligned}
&\mathrel{\int}_0^1\bigg(\dot{f}(t)\ln\bigg(\dfrac{\dot{f}(t)(\lambda+\mu)}{\alpha\lambda}\bigg)
                              -\dot{f}(t)+\dfrac{\alpha\lambda}{\lambda+\mu}\bigg)dt, & &\text{ if}\ f\in \mathbb{AC}_0^M[0,1],\\
                                &\infty, & &\text{ otherwise}.\\
                              \end{aligned}
                              \right.
$$
\end{lemma}

\noindent
\DD.  According to the well-known results (as in
\cite{Mog1}, \cite{BorMog3}, \cite[pp. 13--14]{Feng}),
it is sufficient to show that the Legendre transformation
of the exponential moment of the random variable $\nu_1(1)$ has the following form
$$
\Lambda(x)=\sup\limits_{y\in \mathbb{R}} \bigl(xy-\ln\mathbf{E}e^{y\nu_1(1)} \bigr)=
x\ln\Bigl(\frac{x(\lambda+\mu)}{\alpha\lambda}\Bigr)-x+\frac{\alpha\lambda}{\lambda+\mu}, \ x\geq 0.
$$
Since
$$
\mathbf{E}e^{y\nu_1(1)}=\exp\Bigl\{\frac{\alpha\lambda}{\lambda+\mu} e^{y}-\frac{\alpha\lambda}{\lambda+\mu}\Bigr\},
$$
the applications of the methods of the differential calculus complete
the proof. $\Box$

\begin{lemma} \label{l5.4}
The inequality
\beq \label{5.4}
\mathbf{P}(\nu_2(T)-\nu_2(\delta T)\leq cT)\leq \exp\bigg\{-\frac{\alpha\mu(1-\delta)}{\lambda+\mu}T
 +\frac{\alpha\mu(1-\delta)c}{\lambda+\mu}T-T c\ln c \bigg\}.
\eeq
holds for all $c\in[0,1)$, $\delta\in [0,1]$.
\end{lemma}

\noindent
\DD. Using the Chebyshev inequality, for all $r>0$ we obtain
$$
\begin{aligned}
\mathbf{P} & (\nu_2(T)-\nu_2(\delta T)\leq cT)=
\mathbf{P}(\exp\{-r(\nu_2(T)-\nu_2(\delta T))\}\geq \exp\{-r cT\})
\\ & \leq\frac{\mathbf{E}\exp\{-r(\nu_2(T)-\nu_2(\delta T))\}}{\exp\{-r cT\}}=
\exp\bigg\{e^{-r}\frac{\alpha\mu(1-\delta)}{\lambda+\mu}T-\frac{\alpha\mu(1-\delta)}{\lambda+\mu}T+r cT\bigg\}.
\end{aligned}
$$
Choosing $r=-\ln c$, we obtain the inequality (\ref{5.4}). $\Box$

\begin{lemma} \label{l5.5} The inequality
$$\mathbf{P}\bigg(\bigcap\limits_{k=1}^{n(\varepsilon)-1} B_k \  \bigg| \
\xi_{T}^+(\cdot)\in U_{\varepsilon^3}(g^*),
\bigcap\limits_{k=1}^{n(\varepsilon)} A_k\bigg)\geq \bigg(\frac{\varepsilon^3}{2\Delta(g^*(1)+\varepsilon^3)}\bigg)^{n(\varepsilon)-1},
$$
holds with $g^*$,  $A_k$, $B_k$, and $n(\varepsilon)$ as defined in Section 5.
\end{lemma}

\noindent
\DD.  We show that for $1\leq k \leq n(\varepsilon)-1$, the inequality
\beq \label{5.5}
\mathbf{P}_k:=\mathbf{P}\bigg( B_k \  \bigg| \
\xi_{T}^+(\cdot)\in U_{\varepsilon^3}(g^*),
\bigcap\limits_{k=1}^{n(\varepsilon)} A_k,B_1,\dots,B_{k-1}\bigg)\geq\frac{\varepsilon^3}{2\Delta(g^*(1)+\varepsilon^3)}
\eeq
holds. If the events $\{\omega:\xi_{T}^+(\cdot)\in U_{\varepsilon^3}(g^*)\}$,
$\bigcap\limits_{k=1}^{n(\varepsilon)} A_k$, $B_1,\dots,B_{k-1}$ have all occurred, then

\beq\label{5.6}
\begin{aligned}
&T(g^*(1)+\varepsilon^3)>T(g^*(\tau_{k})+\varepsilon^3)>\xi(\tau_{k}-)\geq T\bigg(\xi_T^+(\tau_{k}-)-(k-1)\bigg(\frac{M}{n(\varepsilon)}-\varepsilon^3\bigg)\bigg)
\\
&>T\bigg(g^*(\tau_{k})-\varepsilon^3-(k-1)\bigg(\frac{M}{n(\varepsilon)}-\varepsilon^3\bigg)\bigg)
>T\bigg(g^*(t_{k})-\varepsilon^3-(k-1)\bigg(\frac{M}{n(\varepsilon)}-\varepsilon^3\bigg)\bigg)
\\
&>T\bigg(g^*(t_{k})-f(t_{k})-\varepsilon^3-(k-1)\bigg(\frac{M}{n(\varepsilon)}-\varepsilon^3\bigg)\bigg)
>\frac{TM}{n(\varepsilon)}-T\varepsilon^3.
\end{aligned}
\eeq
Note that, by definition, the family of
random variables $\zeta_k(m_k)$, $m_k\in \mathbb{N}$
is independent of $\xi^+(t)$ and $\nu_2(t)$, $\zeta_{k-1}(m_{k-1})$, $m_{k-1}\in \mathbb{N}$, $\dots$,
  $\zeta_1(m_1)$, $m_1\in \mathbb{N}$. Hence the random variables $\zeta_k(m_k)$ are independent of
  $\xi(\tau_k-),\dots,\xi(\tau_1-)$.

Therefore, for sufficiently large $T$, using the inequality (\ref{5.6}) and condition $\mathbf{U}$, we obtain
\footnotesize
$$
\begin{aligned}
& \mathbf{P}_k =\mathbf{P}\bigg( B_k \  \bigg| \
\xi_{T}^+(\cdot)\in U_{\varepsilon^3}(g^*),
\bigcap\limits_{k=1}^{n(\varepsilon)} A_k,B_1,\dots,B_{k-1}\bigg)
\\
&=\mathbf{P}\bigg( \zeta_k(\xi(\tau_k-))\in
\bigg(\frac{TM}{n(\varepsilon)}-2T\varepsilon^3;\frac{TM}{n(\varepsilon)}-T\varepsilon^3\bigg) \  \bigg| \
\xi_{T}^+(\cdot)\in U_{\varepsilon^3}(g^*),
\bigcap\limits_{k=1}^{n(\varepsilon)} A_k,B_1,\dots,B_{k-1}\bigg)
\\
&=\sum\limits_{r= \lfloor \frac{TM}{n(\varepsilon)}-T\varepsilon^3 \rfloor }^{\lfloor T(g^*(1)+\varepsilon^3) \rfloor+1}
\mathbf{P}\bigg( \zeta_k(r)\in
\bigg(\frac{TM}{n(\varepsilon)}-2T\varepsilon^3;\frac{TM}{n(\varepsilon)}-T\varepsilon^3\bigg) \  \bigg| \
\xi_{T}^+(\cdot)\in U_{\varepsilon^3}(g^*),
\bigcap\limits_{k=1}^{n(\varepsilon)} A_k, \xi(\tau_k-)=r,B_1,\dots,B_{k-1}\bigg)
\\
& \qquad \qquad \qquad \  \times\mathbf{P}\bigg(\xi(\tau_k-)=r \  \bigg| \
\xi_{T}^+(\cdot)\in U_{\varepsilon^3}(g^*),
\bigcap\limits_{k=1}^{n(\varepsilon)} A_k ,B_1,\dots,B_{k-1}\bigg)
\\
& =\sum\limits_{r= \lfloor \frac{TM}{n(\varepsilon)}-T\varepsilon^3 \rfloor }^{\lfloor T(g^*(1)+\varepsilon^3) \rfloor+1}
\mathbf{P}\bigg( \zeta_k(r)\in
\bigg(\frac{TM}{n(\varepsilon)}-2T\varepsilon^3;\frac{TM}{n(\varepsilon)}-T\varepsilon^3\bigg)\bigg)
\mathbf{P}\bigg(\xi(\tau_k-)=r \  \bigg| \
\xi_{T}^+(\cdot)\in U_{\varepsilon^3}(g^*),
\bigcap\limits_{k=1}^{n(\varepsilon)} A_k,B_1,\dots,B_{k-1}\bigg)
\\
&\geq
\sum\limits_{r= \lfloor \frac{TM}{n(\varepsilon)}-T\varepsilon^3 \rfloor }^{\lfloor T(g^*(1)+\varepsilon^3) \rfloor+1}
\mathbf{P}\bigg( \zeta_k(\lfloor T(g^*(1)+\varepsilon^3) \rfloor+1)\in
\bigg(\frac{TM}{n(\varepsilon)}-2T\varepsilon^3;\frac{TM}{n(\varepsilon)}-T\varepsilon^3\bigg)\bigg)
\\
& \qquad \qquad \qquad\ \times\mathbf{P}\bigg(\xi(\tau_k-)=r \  \bigg| \
\xi_{T}^+(\cdot)\in U_{\varepsilon^3}(g^*),
\bigcap\limits_{k=1}^{n(\varepsilon)} A_k,B_1,\dots,B_{k-1}\bigg)
\\
&=\mathbf{P}\bigg( \zeta_k(\lfloor T(g^*(1)+\varepsilon^3) \rfloor+1)\in
\bigg(\frac{TM}{n(\varepsilon)}-2T\varepsilon^3;\frac{TM}{n(\varepsilon)}-T\varepsilon^3\bigg)\bigg)
\\
&\geq\frac{[T\varepsilon^3]}{\Delta(\lfloor T(g^*(1)+\varepsilon^3) \rfloor+1)}\geq\frac{\varepsilon^3}{2\Delta(g^*(1)+\varepsilon^3)}.
\end{aligned}
$$
\normalsize
Thus, the inequality (\ref{5.5}) is proved. Using the inequality
(\ref{5.5}), we get
$$
\mathbf{P}\bigg(\bigcap\limits_{k=1}^{n(\varepsilon)-1} B_k \  \bigg| \
\xi_{T}^+(\cdot)\in U_{\varepsilon^3}(g^*),
\bigcap\limits_{k=1}^{n(\varepsilon)} A_k\bigg)=
\prod\limits_{k=1}^{n(\varepsilon)-1}\mathbf{P}_k\geq
\bigg(\frac{\varepsilon^3}{2\Delta(g^*(1)+\varepsilon^3)}\bigg)^{n(\varepsilon)-1}.
$$
$\Box$

\begin{lemma} \label{l5.6} The inequality
$$\mathbf{P}\bigg(\bigcap\limits_{l=1}^{\lfloor cT \rfloor}H_{k_l},\bigcap\limits_{l=1}^{\lfloor cT \rfloor}G_{k_l} \
\bigg| \ \nu_2(T)- \nu_2(\delta T)=r\bigg)\leq
\bigg(\frac{\Delta\lfloor 2T\varepsilon \rfloor}{\lfloor T(m_\delta-\varepsilon) \rfloor}\bigg)^{\lfloor cT \rfloor}
$$
holds with $H_{k_l}$,  $G_{k_l}$, $r$, $m_\delta$ defined as in Section 5.
\end{lemma}
\noindent
\DD. Set $G_{k_0}=H_{k_0}=\Omega$.
We will show that the inequality
\beq \label{5.7}
\mathbf{P}_l:=\mathbf{P}\bigg(H_{k_l},G_{k_l} \
\bigg| \ \nu_2(T)- \nu_2(\delta T)=r,\bigcap\limits_{d=0}^{l-1}G_{k_d},\bigcap\limits_{d=0}^{l-1}H_{k_d}\bigg)
\leq \frac{\Delta\lfloor 2T\varepsilon \rfloor}{\lfloor T(m_\delta-\varepsilon) \rfloor}
\eeq
holds for all $1\leq l \leq \lfloor cT \rfloor$.
Now we bound $\mathbf{P}_l$ from above
\footnotesize
$$
\begin{aligned}
\mathbf{P}_l &=\mathbf{P}\bigg(H_{k_l} \
\bigg| \  \nu_2(T)- \nu_2(\delta T)=r, \bigcap\limits_{d=0}^{l}G_{k_d},\bigcap\limits_{d=0}^{l-1}H_{k_d}\bigg)
\mathbf{P}\bigg(G_{k_l} \
\bigg| \ \nu_2(T)- \nu_2(\delta T)=r, \bigcap\limits_{d=0}^{l-1}G_{k_d},\bigcap\limits_{d=0}^{l-1}H_{k_d}\bigg)
\\
&\leq\mathbf{P}\bigg(\zeta_{k_l}(\xi(\tau_{k_l}-))< 2 T\varepsilon \
\bigg| \ \nu_2(T)- \nu_2(\delta T)=r, \bigcap\limits_{d=0}^{l}G_{k_d},\bigcap\limits_{d=0}^{l-1}H_{k_d}\bigg).
\end{aligned}
$$
\normalsize
We note that, by definition, a family of random variables $\zeta_{k_l}(m_{k_l})$, $m_{k_l}\in \mathbb{N}$
is independent of $\xi^+(t)$ and $\nu_2(t)$, $\zeta_{k_{l-1}}(m_{k_{l-1}})$, $m_{k_{l-1}}\in \mathbb{N}$, $\dots$,
  $\zeta_{k_1}(m_{k_1})$, $m_{k_1}\in \mathbb{N}$. Hence $\zeta_{k_l}(m_{k_l})$ are independent of $\xi(\tau_{k_1}-),\dots,\xi(\tau_{k_l}-)$.
Therefore, condition $\mathbf{U}$ implies that  the inequality
\footnotesize
$$
\begin{aligned}
\mathbf{P}&\bigg(\zeta_{k_l}(\xi(\tau_{k_l}-))< 2 T\varepsilon \
\bigg| \ \nu_2(T)- \nu_2(\delta T)=r, \bigcap\limits_{d=0}^{l}G_{k_d},\bigcap\limits_{d=0}^{l-1}H_{k_d}\bigg)
\\
\leq&\sum\limits_{v=\lfloor T(m_\delta-\varepsilon) \rfloor}^{\lfloor T(M+\varepsilon) \rfloor}\mathbf{P}\bigg(\zeta_{k_l}(\xi(\tau_{k_l}-))< 2 T\varepsilon \
\bigg| \ \nu_2(T)- \nu_2(\delta T)=r, \bigcap\limits_{d=0}^{l}G_{k_d},\bigcap\limits_{d=0}^{l-1}H_{k_d},\xi(\tau_{k_l}-)=v\bigg)
\\
&
\times
\mathbf{P}\bigg(\xi(\tau_{k_l}-)=v \
\bigg| \ \nu_2(T)- \nu_2(\delta T)=r, \bigcap\limits_{d=0}^{l}G_{k_d},\bigcap\limits_{d=0}^{l-1}H_{k_d}\bigg)
\\
=&\sum\limits_{v=\lfloor T(m_\delta-\varepsilon) \rfloor}^{\lfloor T(M+\varepsilon) \rfloor}\mathbf{P}\bigg(\zeta_{k_l}(v)< 2 T\varepsilon \bigg)
\mathbf{P}\bigg(\xi(\tau_{k_l}-)=v \
\bigg| \ \nu_2(T)- \nu_2(\delta T)=r, \bigcap\limits_{d=0}^{l}G_{k_d},\bigcap\limits_{d=0}^{l-1}H_{k_d}\bigg)
\\
\leq& \sum\limits_{v=\lfloor T(m_\delta-\varepsilon) \rfloor}^{\lfloor T(M+\varepsilon) \rfloor}\mathbf{P}\bigg(\zeta_{k_l}(\lfloor T(m_\delta-\varepsilon) \rfloor)< 2 T\varepsilon\bigg)\mathbf{P}\bigg(\xi(\tau_{k_l}-)=v \
\bigg| \ \nu_2(T)- \nu_2(\delta T)=r, \bigcap\limits_{d=0}^{l}G_{k_d},\bigcap\limits_{d=0}^{l-1}H_{k_d}\bigg)
\\
=&\mathbf{P}\bigg(\zeta_{k_l}(\lfloor T(m_\delta-\varepsilon) \rfloor)< 2 T\varepsilon\bigg)
\leq \frac{\Delta[2 T\varepsilon]}{\lfloor T(m_\delta-\varepsilon) \rfloor}
\end{aligned}
$$
\normalsize
 for sufficiently small $\varepsilon$.
Thus, the inequality (\ref{5.7}) is proved. Using the inequality
(\ref{5.7}), we obtain
$$
\mathbf{P}\bigg(\bigcap\limits_{l=1}^{\lfloor cT \rfloor}H_{k_l},\bigcap\limits_{l=1}^{\lfloor cT \rfloor}G_{k_l} \
\bigg| \ \nu_2(T)- \nu_2(\delta T)=r\bigg)=
\prod\limits_{l=1}^{\lfloor cT \rfloor}\mathbf{P}_l\leq
\bigg(\frac{\Delta[2 T\varepsilon]}{\lfloor T(m_\delta-\varepsilon) \rfloor}\bigg)^{\lfloor cT \rfloor}.
$$
$\Box$

\section*{Acknowledgments}
This work is supported by FAPESP grant 2017/20482-0.

AL supported by RSF project 18-11-00129. AL thanks the Institute of Mathematics and
Statistics of University of S\~ao Paulo for hospitality. AY thanks CNPq and
FAPESP for the financial support via the grants 301050/2016-3 and 2017/10555-0, respectively.

\end{document}